\documentclass[letterpaper, 10pt, conference]{ieeeconf}
\IEEEoverridecommandlockouts               
\usepackage{graphicx}
\usepackage{epsfig,epstopdf}
\graphicspath{{../pdf/}{../jpeg/}}
\DeclareGraphicsExtensions{.pdf,.jpeg,.png}
\usepackage{amsmath,amssymb,amsfonts,mathrsfs}
\usepackage{multirow}
\usepackage{booktabs}
\usepackage{color}
\usepackage{algpseudocode}
\usepackage{tikz}
\usepackage{color}
\usepackage{cite}
\usepackage{setspace}
\usepackage{algorithm}  
\usepackage{hyperref}
\hypersetup{
colorlinks=true,
linkcolor=black,
filecolor=magenta,      
urlcolor=cyan,
pdftitle={Overleaf Example},
pdfpagemode=FullScreen,
}

\makeatletter
\def\algbackskip{\hskip-\ALG@thistlm}
\makeatother

\allowdisplaybreaks[4]
\usepackage{theorem}
\newtheorem{theorem}{Theorem}

\newtheorem{lemma}{Lemma}

\newtheorem{property}{Property}
\newtheorem{definition}{Definition}

\newtheorem{assumption}{Assumption}
\newtheorem{remark}{Remark}
\newtheorem{example}{Example}

\DeclareMathOperator*{\argmin}{arg\,min}

\title{\bf 
Safe Zeroth-Order Convex Optimization Using Quadratic Local Approximations
}

\author{
\thanks{This work was supported by the Swiss National Science Foundation under the NCCR Automation (grant agreement 51NF40\textunderscore 80545).}
\thanks{Baiwei Guo and Giancarlo Ferrari-Trecate are with the DECODE group, Institute of Mechanical Engineering, EPFL, Switzerland. email: {\tt \{baiwei.guo, giancarlo.ferraritrecate\}@epfl.ch}.} 
\thanks{Yuning Jiang is with the Predict group, IGM, EPFL, Switzerland. email: {\tt \{yuning.jiang\}@ieee.org}.}
Baiwei Guo, Yuning Jiang, Maryam Kamgarpour, Giancarlo Ferrari-Trecate
\thanks{Maryam Kamgarpour is with the Systems Control And Multiagent Optimization Research (sycamore) labs at IGM, EPFL, email:{\tt \{maryam.kamgarpour\}@epfl.ch}.%
}}

\begin{document}
\maketitle

\begin{abstract}
We address black-box convex optimization problems, where the objective and constraint functions are not explicitly known but can be sampled within the feasible set. The challenge is thus to generate a sequence of feasible points converging towards an optimal solution. By leveraging the knowledge of the smoothness properties of the objective and constraint functions, we propose a novel zeroth-order method, SZO-QQ, that iteratively computes quadratic approximations of the constraint functions, constructs local feasible sets and optimizes over them. We prove convergence of the sequence of the objective values generated at each iteration to the minimum. Through experiments, we show that our method can achieve faster convergence compared with state-of-the-art zeroth-order approaches to convex optimization.
\end{abstract}

\section{Introduction}
\label{sec: introduction}
Applications ranging from energy resources operation \cite{chu2021frequency}, machine learning \cite{chen2020DFO} and trajectory optimization \cite{manchester2016derivative} to optimal control tuning~\cite{xu2021vabo} often require solving complex optimization problems where the fulfillment of the hard constraints is essential. In case the modelling of the objective and constraint functions is not available or the given data from safe system operation is not sufficient for modelling, zeroth-order optimization methods can be used to derive the optimum without violation of constraints. The methods of this category only rely on the possibility of sampling the unknown objective and constraint functions evaluated at a set of chosen points~\cite{Bajaj2021}. The decision variables are then iteratively updated using the collected evaluations. 

Classical techniques for zeroth-order unconstrained optimization can be classified as direct-search-based (where a set of points around the current point is searched for a lower value of the objective function) and model-based (where a local model of the objective function around the current point is built and used for local optimization)~\cite[Chapter~9]{jorge2006numerical}. Examples for these two categories are, respectively, pattern search methods \cite{lewis2000pattern} and trust region methods \cite{conn2000trust}. For constrained problems, a common approach is to include in the cost barrier functions penalizing constraint violations so as to cast the constrained optimization problem into an unconstrained one \cite{lewis2002globally,audet2009progressive}. By adopting log barrier functions to penalize the proximity towards the boundary of the feasible set, Extremum Seeking methods estimate the gradient of the new objective function by adding sinusoidal perturbations to the decision variables~\cite{arnold2015model}. However, most penalty-based methods lack a mechanism to tune the penalty coefficients for ensuring feasible sampling. 

To guarantee feasible sampling with a high probability,~\cite{sui2015safe} proposes an iterative algorithm, called SafeOpt, that expands the estimated feasible set by modeling the constraints as Gaussian processes. However, the computational complexity of this algorithm scales exponentially with the number of the decision variables. To address this issue, a gradient-descent penalty-based approach, called s0-LBM, is proposed in~\cite{usmanova2019log}. This method exploits the constraint functions' smoothness to calculate the step size of the descent. Although this method can address non-convex problems and comes with a worst-case complexity that is polynomial in problem dimension, it might converge slowly, {\color{black}even for convex problems}. The reason is that as the iterates approach the boundary of the feasible set, the log-barrier function and its derivative become very large, leading to very conservative local feasible sets that result in slow progress of the iterates. The work \cite{usmanova2019safe} addresses the zeroth-order optimization of convex objectives subject to unknown polytopic constraints. It ensures sampling feasibility by estimating the polytope iteratively and constructing confidence bounds for the estimates. While the optimization algorithm in \cite{usmanova2019safe}  enjoys faster convergence than \cite{usmanova2019log}, the approach is not readily generalizable to  nonlinear constraints.  

\textit{Contributions:} In this work, we aim to solve convex optimization problems based on Lipschitz and smoothness constants instead of the exact knowledge of the objective and constraint functions. To achieve fast convergence, we approximate the constraint functions with known quadratic functions (in contrast to the affine approximation used in \cite{sui2015safe,usmanova2019log}) and tackle the optimization problem by solving a series of known Quadratically Constrained Quadratic Programs (QCQP). Our main contributions are summarized as follows:
\begin{enumerate}
    \item a method to construct local feasible sets, which is {less conservative than those in \cite{sui2015safe,usmanova2019log} when the subproblem solutions are close to the boundary of the original feasible set}. To this purpose we use local quadratic approximations of the constraints based on the knowledge of the Lipschitz constants of the constraints' gradients;
    \item a zeroth-order optimization algorithm where sampling feasibility is guaranteed and the objective function values provably converge to the minimum;
    \item numerical examples comparing the speed of convergence of our method with state-of-the-art alternative approaches as well as an application of the proposed algorithm to an optimal control problem.
\end{enumerate}

\textit{Notations:}
We use $e_i\in\mathbb R^d$ to define the $i$-th standard basis of vector space $\mathbb{R}^d$. Given a vector $x\in \mathbb{R}^d$ and a scalar $\epsilon>0$, we write $x = [x^{(1)}, \ldots, x^{(d)}]^\top$ and $\mathcal B_\epsilon(x)=\{y:\|y-x\|\leq \epsilon\}$, where $\|\cdot\|$ is the 2-norm. We define $\mathbb{Z}_i^j=\{i,i+1,\ldots,j\}$.

\section{Problem Formulation}
\label{sec::Prob}
We address the constrained convex optimization problem
\begin{equation}
\label{eq: optimization problem}
 \min_{x\in \mathbb{R}^d}  \;\; f_0(x)\quad \text{subject to}\;\;
x\in \Omega,
\end{equation}
where $\Omega:=\{x\in \mathbb{R}^d: f_i(x)\leq 0, i\in\mathbb Z_1^m\}$ is the feasible set.
We consider the setting where the convex functions $f_i:\mathbb R^d\to\mathbb R$, $i\in \mathbb{Z}^n_1$, are not explicitly known but they can be queried  within $\Omega$. We further assume, without loss of generality, the convex objective function $f_0(x)$ is explicitly known and linear. Indeed, when the convex function $f_0(x)$ in \eqref{eq: optimization problem} is not known but can be queried, the problem in \eqref{eq: optimization problem} can be written as 
\begin{equation*}
\begin{aligned}
 \min_{(x,\gamma)\in \mathbb{R}^{d+1}
} \quad & \gamma\\
\mathrm{subject\;to}  \quad & f_0(x)-\gamma\leq 0,\\
 & f_i(x)\leq 0, \quad i \in\mathbb Z_1^m,
\end{aligned}
\end{equation*}
where the objective function is now known and linear. 

We also consider the following assumptions on the smoothness of the objective and constraint functions and the availability of a strictly feasible point.

\begin{assumption}
\label{ass: smoothness}
The functions $f_i(x)$, $i\in\mathbb Z_0^m$ are continuously differentiable and satisfy for any $x_1$, $x_2\in \mathbb{R}^d$,
\begin{equation}
    \begin{aligned}
    |f_i(x_1)-f_i(x_2)|&\leq L_i\|x_1-x_2\|,\\
    \|\nabla f_i(x_1)-\nabla f_i(x_2)\|&\leq M_i\|x_1-x_2\|.\\
    \end{aligned}
    \label{eq: smoothness}
\end{equation}
\end{assumption}
In the following, we use the notation $L_{\max}: = \max_{i\geq 1}L_i$, $M_{\max}: = \max_{i\geq 1}M_i$, $L_{i,\text{inf}}:=\inf\{L:\| f_i(x_1)- f_i(x_2)\| \leq L\|x_1-x_2\|,\forall x_1,x_2\in \Omega\}$ and $M_{i,\text{inf}}:=\inf\{M:\|\nabla f_i(x_1)-\nabla f_i(x_2)\| \leq M\|x_1-x_2\|,\forall x_1,x_2\in \Omega\}$.
\begin{remark}
The bounds in \eqref{eq: smoothness} are utilized in several works on zeroth-order optimization, e.g., \cite{wibisono2012finite,conn2009introduction}. As it will be clear in the sequel, these bounds allow one to estimate the error of local approximations of the unknown functions and their derivatives. Moreover, our method does not require to precisely know $L_{i,\text{inf}}$ and $M_{i,\text{inf}}$, but only upper bounds to their values, see Remark \ref{rmk: LM}.

\end{remark}
\begin{assumption}
\label{ass: strict feasible}
There exists a known strictly feasible point $x_0$, i.e., $f_i(x_0)<0$ for all $i\in\mathbb Z_1^m$. 
\end{assumption}
\begin{remark}
The existence of a strictly feasible point is called Slater's Condition and commonly assumed in several optimization methods \cite{boyd2004convex}. Moreover, several works on safe learning \cite{sui2015safe,usmanova2019safe} assume that a strictly feasible point is given to initialize the algorithm. This is necessary for designing an algorithm whose iterates remain feasible since the constraint functions are unknown a priori. Practically, this assumption holds in several applications. For example, in any robot mission planning,  the robot is placed initially at a safe point and needs to gradually explore the neighboring regions while ensuring feasibility of its trajectory. Similarly, in the optimization of manufacturing processes, often an initial set of (suboptimal) design parameters that satisfy the problem constraints are known \cite{rupenyan2021performance}.
\end{remark}
\begin{assumption}
\label{ass: bounded_sublevel}
There exists $\beta\in \mathbb{R}$ such that the sublevel set $\mathcal{P}_\beta=\{x\in \Omega:f_0(x)\leq \beta\}$ is bounded and includes the initial feasible point $x_0$.
\end{assumption}
The above assumption ensures the iterates remain bounded and have a convergent subsequence (due to compactness of $\mathcal{P}_\beta$) as will be shown in Theorem \ref{thm: main_theorem}. 

Based on the aforementioned assumptions and properties, in the next section we design an algorithm that iteratively optimizes $f_0(x)$ over local safe sets constructed by exploiting Assumption  \ref{ass: smoothness} on smoothness.

\section{The Proposed Zeroth-Order Algorithm}
\label{sec: algorithm_tot}
Here, we propose an approach to construct  local feasible sets around a given feasible point using function queries. To do so, we first recall properties of a gradient estimator constructed through finite differences.

The gradients of the unknown functions $\{f_i\}_{i=1}^m$ can be approximated through finite differences as  
\begin{equation}
\label{eq: gradient approximation}
{\nabla}^\nu f_i\left(x\right):=\sum_{j=1}^{d} \frac{f_i\left(x+\nu e_{j}\right)-f_i\left(x\right)}{\nu} e_{j}
\end{equation}
with $\nu>0$. The estimation error is denoted as $$\Delta^\nu_i(x):={\nabla}^\nu f_i\left(x\right)-\nabla f_i(x).$$
From \eqref{eq: smoothness}, we have the following result.
\begin{lemma}[\cite{berahas2021theoretical}, Theorem 3.2]
\label{lmm: gradient estimation error}
Under Assumption \ref{ass: smoothness}, the gradient estimation error at $x\in\mathbb{R}^d$ can be bounded as
\begin{equation}
\left\|\Delta^\nu_i(x)\right\|_{2} \leq \alpha_i \nu, \text{ with } \alpha_i = \frac{\sqrt{d}  M_i}{2}.
\label{eq: gradient approximation error}
\end{equation}
\end{lemma}

\subsection{Local feasible set construction}
\label{sec: Safesetconstruction}
Based on the estimated gradient and the error bound above, we build a local feasible set around a feasible point $x_0$. 
\begin{theorem}
\label{thm: safesetconstruction}
For any strictly feasible point $x_0$, let $l^* = \min_{i=1, \dots, m}-f_i(x_0)/L_{\mathrm{max}} \text{ and } \nu^* = 2l^*/\sqrt{d}$. Define 
\begin{equation}
\label{eq:safe_set}
\begin{aligned}
\mathcal{S}^\circ_i(x_0) :=&\big\{x: f_i(x_0)+{\nabla}^{\nu^*} f_i\left(x_0\right)^\top(x-x_0)+\\
&\qquad\qquad\qquad\quad 2M_i\|x-x_0\|^2\leq 0 \big\}.  
\end{aligned}
\end{equation}
Under Assumption \ref{ass: smoothness}, the set $\mathcal{S}^\circ(x_0):=\bigcap^m_{i=1}\mathcal{S}^\circ_i(x_0)$  satisfies $\mathcal{S}^\circ(x_0) \subset \Omega$.
\end{theorem}
\textbf{Proof.}
We first partition $\mathcal{S}^\circ_i(x_0)$ as 
\[
\mathcal{S}^\circ_i(x_0) = \left(\mathcal{S}^\circ_i(x_0)\bigcap \mathcal B_{l^*}(x_0)\right)\bigcup\left(\mathcal{S}^\circ_i(x_0)\setminus \mathcal B_{l^*}(x_0)\right)
\] 
and notice that $\mathcal{S}^\circ_i(x_0)\bigcap \mathcal B_{l^*}(x_0)\subseteq \Omega$. Then, it only remains to check that $\mathcal{S}^\circ_i(x_0)\setminus \mathcal B_{l^*}(x_0)\subseteq \Omega.$

For $x\in \mathcal{S}^\circ_i(x_0)\setminus B_{l^*}(x_0)$, we have $\frac{\sqrt{d}\nu^*}{2}= l^* \leq\|x-x_0\|.$
By the mean value theorem, for any $i$ there exists $\theta_i\in [0,1]$ such that \begin{equation}
\label{eq: meanvaluetheorem}
\begin{aligned}
    \hspace{-0.3mm}f_i(x)
    =&  f_i(x_0)+\nabla f_i\left(x_0 + \theta_i(x-x_0)\right)^\top(x-x_0)\\
    = &f_i(x_0) + \nabla f_i\left(x_0 \right)^\top(x-x_0) +\\
    & \left(\nabla f_i\left(x_0 + \theta_i(x-x_0)\right) - \nabla f_i\left(x_0 \right)\right)^\top(x-x_0)\\
    \leq& f_i(x_0) + {\nabla}^{\nu^*} f_i\left(x_0\right)^\top(x-x_0)\\
    &+ \frac{\sqrt{d}\nu^* M_i}{2}\|x-x_0\| + M_i\|x-x_0\|^2\\
    \leq& f_i(x_0) + {\nabla}^{\nu^*} f_i\left(x_0\right)^\top(x-x_0) + 2M_i\|x-x_0\|^2\\
    \leq& 0,
\end{aligned}
\end{equation}
where the first inequality is due to Assumption \ref{ass: smoothness} while the second can be derived according to the definition of $\nu^*$ in Theorem \ref{thm: safesetconstruction}. Hence, $\mathcal{S}^\circ_i(x_0)\setminus \mathcal B_{l^*}(x_0)\subseteq \Omega.$ Since $\mathcal{S}^\circ_i(x_0)\subseteq \Omega, \forall i$, then $\mathcal{S}^\circ(x_0)\subseteq \Omega$. \hfill$\blacksquare$

By construction, we see that if $x_0$ is strictly feasible, then it belongs to $\mathcal{S}^\circ(x_0)$. Moreover, the set $\mathcal{S}^\circ(x_0)$ is convex since $\mathcal{S}^\circ(x_0) = \cap^{m}_{i=1}\mathcal{S}^\circ_i(x_0)$ and, for any $i$, $\mathcal{S}^\circ_i(x_0)$ is a $d$-dimensional ball. In the sequel, we call $\mathcal{S}^\circ(x_0)$ a \textit{local feasible set around} $x_0$.

\begin{remark}
\label{rmk: remark on the local estimation}
The works ~\cite{usmanova2019log,sui2015safe} consider an alternative approximation of the functions and a local feasible set 
\[
\mathcal{C}(x_0):=\bigcap^m_{i=1}\left\{x:\|x-x_0\|\leq -\frac{f_i(x_0)}{L_{\max}}\right\}.
\]
We see that $\mathcal{C}(x_0)= \{x:F^L_i(x)\leq0,\forall i \}$ where $F^L_i(x):=f_i(x_0)+L_{\max}\|x-x_0\|$. In contrast, $\mathcal{S}^\circ(x_0)=\{x:F^M_i(x)\leq0,\forall i \}$ where $F^M_i(x):=f_i(x_0)+{\nabla}^{\nu^*} f_i\left(x_0\right)^\top(x-x_0)+2M_{i}\|x-x_0\|^2$. One can easily verify that $F^M_i(x)<F^L_i(x)$ when $\|x-x_0\|$ is sufficiently small. If $x_0$ is close to the boundary of the feasible set, the distance $\|x-x_0\|$ is small for any $x\in \mathcal{C}(x_0)$ and thus $\mathcal{C}(x_0)\subset \mathcal{S}^\circ(x_0)$ due to $F^M_i(x)<F^L_i(x)$. We illustrate this point in Example~\ref{ex: lessconservativeness}. The above analysis reveals why $\mathcal{S}^\circ(x_0)$ can be less conservative than $\mathcal{C}(x_0)$, which constitutes the main reason for which our method can achieve faster convergence than \cite{usmanova2019log}.\end{remark}
\begin{example}
\label{ex: lessconservativeness}
Consider the optimization problem \eqref{eq: optimization problem} with $d=1$, $m=2$, $f_0(x) = 3x$, $f_1(x) = x^2-x-0.75$ and $x_0 = 1.49$, which is close to the boundary of the feasible set $[-1,1.5]$. We let ${L_{\max}}=3.01$ and ${M_{\max}} = 3$, then $\mathcal{T}(x_0)=[1.4836,1.4964] \subset \mathcal{S}^\circ(x_0) = [1.1502,1.4998]$.
\end{example}

\subsection{The proposed algorithm}
\label{sec: algorithm}
The proposed method to solve problem \eqref{eq: optimization problem} is summarized in Algorithm~\ref{alg:sslo}. 
The main idea is to start from the strictly feasible initial point $x_0$ and iteratively solve $\hyperlink{starstar}{(*)}$ in Algorithm~\ref{alg:sslo}, which is a convex Quadratically Constrained
Quadratic Program (QCQP), until the distance $\|x_{k+1}-x_k\|$ is smaller than a given threshold $\xi>0$ (see Line 6-8). By using the penalty term $\mu \|x-x_k\|^2$ in $\hyperlink{starstar}{(*)}$, we can guarantee that the increment $\|x_{k+1}-x_k\|$ diminishes and Algorithm \ref{alg:sslo} terminates as shown later in Theorem \ref{thm: main_theorem}.
\begin{algorithm}[htbp!]
\caption{Safe Zeroth-Order Sequential QCQP (SZO-QQ)}
\label{alg:sslo}
\textbf{Input:} $\mu>0$, $\xi>0$, initial feasible point $x_0\in\Omega$\\
\textbf{Output:} $\tilde{x}$
\begin{algorithmic}[1]
\State Choose $M_i >M_{i,\text{inf}}$, for $i\in \mathbb{Z}^m_1$
\State $k\gets 0, \textsc{TER}= 0$ 
\While{$\textsc{TER} = 0$}
\State Compute $\mathcal{S}(x_k)$ based on \eqref{eq:safe_set} and \eqref{eq: stepsize_2}.
\State 
$x_{k+1} = \argmin_{x\in \mathcal{S}(x_{k})} f_0(x)+\mu \|x-x_k\|^2 \quad \quad \hypertarget{starstar}{(*)} $ 
\If{$\|x_{k+1}-x_{k}\|\leq \xi$}
\State $\tilde{x}\gets x_{k+1}$, $\textsc{TER}\gets 1$
\EndIf
\State $k\gets k+1$
\EndWhile
\end{algorithmic}
\end{algorithm}

The construction of $\mathcal{S}(x_k)$ is based on Theorem \ref{thm: safesetconstruction}. We let $\mathcal{S}(x_k) = \bigcap^{m}_{i = 1}\mathcal{S}_i(x_k)$ where $M_i>M_{i,\mathrm{inf}}$,  
\begin{equation}
  l_k^*= \min_i\frac{-f_i(x_k)}{L_{\max}},\quad \nu_k^* =\min\{ \frac{2l_k^*}{\sqrt{d}},\frac{1}{k}\}.
   \label{eq: stepsize_2}
\end{equation}
and
\begin{equation}
\label{eq:safe_set_alg}
\begin{aligned}
\mathcal{S}_i(x_k) :=&\big\{x: f_i(x_0)+{\nabla}^{\nu_k^*} f_i\left(x_0\right)^\top(x-x_0)+\\
&\qquad\qquad\qquad\quad 2M_i\|x-x_0\|^2\leq 0 \big\}.  
\end{aligned}
\end{equation}
By requiring $M_i>M_{i,\mathrm{inf}}$, we guarantee that $\forall k$, $x_k$ is strictly feasible (see the proof of Theorem \ref{thm: safesetconstruction}). Thanks to \eqref{eq: stepsize_2}, we can ensure that the gradient approximation error goes to 0, which is essential for the convergence analysis in Theorem~\ref{thm: main_theorem}. 

\begin{remark}
\label{rmk: LM}
To run Algorithm \ref{alg:sslo}, we need to set $L_i>L_{i,\mathrm{inf}}$ and $M_i>M_{i,\mathrm{inf}}$ for $i\geq 1$. Here, we discuss the effect of under/overestimating $L_{i,\mathrm{inf}}$ and $M_{i,\mathrm{inf}}$. If $L_i<L_{i,\mathrm{inf}}$ or $M_i<M_{i,\mathrm{inf}}$, the set $\mathcal{S}(x_k)$ built according to \eqref{eq: stepsize_2} and \eqref{eq:safe_set_alg} might not be feasible. Consequently, when infeasible samples happen, one can increase $L_i$ and $M_i$ for $i\geq 1$ by a certain factor. At the same time, one should avoid too large a value for $M_i$, since if $M_i\gg M_{i,\mathrm{inf}}$, the approximation used to construct $\mathcal{S}(x_i)$ can be very conservative and thus the convergence of Algorithm \ref{alg:sslo} can become much slower.  
\end{remark}

Algorithm \ref{alg:sslo} is similar to Sequential QCQP (SQCQP)~\cite{solodov2004sequential}, where at each iteration quadratic proxies for both the objective and constraint functions are built based on the local gradient vectors and Hessian matrices. Applications of SQCQP to optimal control have received increasing attention \cite{messerer2020determining,messerer2021survey}, due to the development of efficient solvers for QCQP subproblems \cite{frison2022introducing}. Different from SQCQP \cite{solodov2004sequential}, Algorithm~\ref{alg:sslo}, which we also call SZO-QQ in the remainders, can guarantee feasible sampling and does not require line search, which is costly due to the required generation of extra examples. Due to the lack of line search, Algorithm \ref{alg:sslo}'s convergence cannot be proved by simply following the arguments in \cite{solodov2004sequential}. 

\section{Theoretical Convergence Analysis}
In this section, we first consider the sequence $\{x_k\}_{k\geq 1}$ generated by Algorithm~\ref{alg:sslo} without the termination condition. We show that under mild assumptions, the sequence $\{f_0(x_k)\}_{k\geq 1}$ converges to the minimum objective value of \eqref{eq: optimization problem}. Then, we discuss the choice of the termination condition and the required number of iterations for the algorithm to stop.
\label{sec: accumulation_point}

\subsection{Objective value converging to minimum}
We let $\xi = 0$, which is to say the termination condition $\|x_{k+1}-x_{k}\|\leq \xi$ in Algorithm \ref{alg:sslo} is omitted. We characterize the limit of the sequence $\{f_0(x_k)\}_{k\geq 1}$ generated by the algorithm in the following theorem.
\begin{theorem}
\label{thm: main_theorem}
Under Assumption \ref{ass: smoothness}-\ref{ass: bounded_sublevel}, if $\xi = 0$, the sequence $\{x_k\}_{k\geq 1}$ in Algorithm~\ref{alg:sslo} has the following properties:
\begin{itemize}
    \item[1.] the sequence $\{f_0(x_k)\}_{k\geq 1}$ is non-increasing;
    \item[2.] the sequence $\{x_k\}_{k\geq 1}$ has at least one accumulation point $x_c$\footnote{Theorem \ref{thm: main_theorem} does not imply that $\{x_k\}_{k\geq 1}$ is convergent. For this property to hold, a sufficient condition is $\sum^{\infty}_{k=0}\|x_{k}-x_{k+1}\|<\infty$. However, in the proof of Theorem \ref{thm: main_theorem}, we have only been able to show a weaker condition, namely, $\sum^{\infty}_{k=0}\|x_{k}-x_{k+1}\|^2<\infty$.};
    \item [3.] for any accumulation point $x_c$ of the sequence, $\lim_{k\rightarrow \infty}f_0(x_k)=f_0(x_c)> -\infty$. 
\end{itemize}
\end{theorem}
The proof of Theorem \ref{thm: main_theorem}, provided in Appendix \ref{sec: 1st argument of thm2}, leverages the optimality of the solution to the subproblem $\hyperlink{starstar}{(*)}$ in Algorithm \ref{alg:sslo}, the boundedness of the sublevel set $\mathcal{P}_\beta$ required by Assumption \ref{ass: bounded_sublevel} and the continuity of $f_0(x)$. 

In view of point 3 of Theorem \ref{thm: main_theorem}, to show that $\{f_0(x_k)\}_{k\geq 1}$ converges to the minimum objective value of \eqref{eq: optimization problem}, we need to show there exists an accumulation point $x_c$ of $\{x_k\}_{k\geq 1}$ that is the optimum of  \eqref{eq: optimization problem}. We now address the optimality of any accumulation point $x_c$ in Theorem \ref{thm: optimality} below.
\begin{assumption}
\label{ass: constraint qualification}
There exists an accumulation point $x_c$ such that $\nabla f_i(x_c)\neq 0,\forall i\in \mathcal{A}(x_c)$, where $\mathcal{A}(x_c):=\{i:f_i(x_c)=0\}$
\end{assumption}
\begin{remark}
For proving the convergence of Interior Point Methods \cite[Theorem 19.1]{jorge2006numerical}, one often assumes the Linear Independence Constraint Qualification (LICQ), requiring that vectors $\nabla f_i(x_c)\text{ for } i\in \mathcal{A}(x_c)$ are linearly independent. Note that Assumption \ref{ass: constraint qualification} is weaker than LICQ because if there exists $i_0\in\mathcal{A}(x_c)$ such that $\nabla f_{i_0}(x_c) = 0$ then the vectors $\nabla f_i(x_c)\text{ for } i\in \mathcal{A}(x_c)$ cannot be linearly independent.
\end{remark}
\begin{theorem}
\label{thm: optimality}
Under Assumptions \ref{ass: smoothness}-\ref{ass: constraint qualification}, any accumulation point $x_c$ is an optimizer of problem \eqref{eq: optimization problem}.
\end{theorem}

The complete proof of Theorem \ref{thm: optimality} can be found in Appendix \ref{sec: proof_of_optimality}. Below, we provide an outline of it, which is based on the  following optimality conditions for the constrained optimization problem \eqref{eq: optimization problem} as stated below.

\begin{definition}[\cite{dutta2013approximate}]
\label{def: KKT}
A pair $(x,\lambda)$ with $x\in \Omega$ and $\lambda\in \mathbb{R}^m_{\geq 0}$ is a $\eta-$approximate KKT (also written as $\eta-$KKT) point with $\eta>0$ of Problem \eqref{eq: optimization problem} if it fulfills
\begin{equation}
\label{eq: kkt}
\begin{aligned}
    \|\nabla f_0(x) + \sum^m_{i=1}\lambda^{(i)}\nabla f_i(x)\| &\leq \eta,\\
    |\lambda^{(i)}f_i(x)| &\leq \eta,\quad i\in \mathbb{Z}^{m}_{1}.\\
\end{aligned}
\end{equation}
If $(x^*,\lambda^*)$ with $x^*\in \Omega$ and $\lambda^*\in \mathbb{R}^m_{\geq 0}$ fulfills \eqref{eq: kkt} with $\eta = 0$, we say that it is a KKT point.
\end{definition}
We notice that due to Slater's Condition, $x^*$ is a minimizer of \eqref{eq: optimization problem} if and only if there exists $\lambda^* \in \mathbb{R}^d_{\geq 0}$ such that $(x^*,\lambda^*)$ is a KKT point of \eqref{eq: optimization problem}.

\vspace{0.3cm}
\noindent\textit{Sketch of the proof of Theorem \ref{thm: optimality}: } Let $x_c$ be an accumulation point of $\{x_k\}_{k\geq 1}$ such that $\nabla f_i(x_c)\neq 0,\forall i\in \mathcal{A}(x_c)$ (see Assumption \ref{ass: constraint qualification}). Since $x_k$ is feasible for any $k$ and $\Omega$, the point $x_c$ belongs to $\Omega$. If $\nabla f_0(x_c) = 0$, then with $\lambda_c = 0$, the pair $(x_c,\lambda_c)$ is a KKT point and thus $x_c$ is a minimizer. Now, we consider $\nabla f_0(x_c) \neq 0$. It follows that that the active set $\mathcal{A}(x_c)=\{i: f_i(x_c)=0\}$ is not empty (see Lemma \ref{prop: active set non-empty} in Appendix \ref{sec: state_proof_lemma_nonempty_active_set} for the rigorous statement). We can then show the optimality of $x_c$ through contradiction as follows. 

We assume $x_c$ is not a minimizer and let $x_s\in \Omega$ be such that $f_0(x_s)<f_0(x_c)$. We denote the
line segment between $x_c$ and $x_s$ as $x(t):=x_c+t(x_s-x_c), t\in [0,1]$. By utilizing Lemma \ref{prop: active set non-empty} in Appendix \ref{sec: state_proof_lemma_nonempty_active_set}, we can find $\bar{t}\in (0,1)$ such that any $0<t_p<\bar{t}$ verifies $x(t_p)$ belonging to the interior of $\mathcal{S}(x_c)$ and $f_0(x(t_p))+\mu \|x(t_p)-x_{c}\|^2<f_0(x_{c})$. For this $t_p$ and any $\tau>0$, there exists $k_p>0$ verifying $\|x_{k_p}-x_c\|\leq\tau$ and $ x(t_p)\in \mathcal{S}(x_{k_p})$ (see Lemma \ref{claim: more_optimal_point} in Appendix \ref{sec: state_proof_more_optimal_point} for the rigorous statement). By letting $\tau$ be sufficiently small, we have $f_0(x(t_p))+\mu \|x(t_p)-x_{k_p}\|^2<f_0(x_{c})$. Since $x_{k_p+1} = \argmin_{x\in \mathcal{S}(x_{k_p})}f_0(x)+\mu \|x-x_{k_p}\|^2$, we have $f_0(x_{k_p+1})+\mu \|x_{k_p+1}-x_{k_p}\|^2\leq f_0(x(t_p))+\mu \|x(t_p)-x_{k_p}\|^2<f_0(x_c)$, which contradicts the monotonicity shown in point 1 of Theorem \ref{thm: main_theorem}. Thus, $x_c$ is a minimizer of \eqref{eq: optimization problem}.

\subsection{Computational complexity analysis}
The result in Theorem \ref{thm: optimality} is asymptotic, but in practice only finitely many iterations can be computed. Algorithm \ref{alg:sslo} adopts the termination condition  $\|x_{k+1}-x_k\|<\xi$ for $\xi > 0$. The motivation for this particular stopping criteria is as follows. With this condition, we expect that the KKT point of the inner optimization problem $\hyperlink{starstar}{(*)}$ in Algorithm \ref{alg:sslo} fulfills the approximate KKT conditions of the original problem, namely \eqref{eq: kkt}. To see this point, for $k=K-1$ we consider the KKT point of the subproblem $(*)$, $(x_K,\lambda_K)$.  If the error term $\|x_K-x_{K-1}\|$ is small, one can show that the left-hand-side terms in \eqref{eq: kkt} evaluated at $(x,\lambda) =(x_K,\lambda_K)$ are also small.

Our next goal is to address whether the termination condition $\|x_{k+1}-x_{k}\|<\xi$ is triggered in a finite number of iterations. And if yes, how many iterations are needed.
\begin{theorem}
\label{thm: termination_complexity}
Under Assumption \ref{ass: smoothness}-\ref{ass: constraint qualification}, the increment $\|x_{k+1}-x_{k}\|$ converges to 0 as $k \rightarrow \infty$. Moreover, $K :=\inf\{k: \|x_{k+1}-x_k\|\leq \xi\}$ satisfies
\begin{equation}
\label{eq: complexity}
    K\leq\overline{K} = \frac{f_0(x_0) - \inf\{f_0(x): x\in \Omega\}}{\mu\xi^2}.
\end{equation}
Hence, Algorithm \ref{alg:sslo} terminates within $\overline{K}$ iterations.
\end{theorem}
\textbf{Proof.}
Since $x_{k+1}= \argmin f_0(x_k)+\mu\|x-x_k\|^2$, we know that $f_0(x_{k})-f_0(x_{k+1})\geq \mu\|x_{k}-x_{k+1}\|^2$. Then, we have
\begin{equation}
\label{eq:square_convergence}
\begin{aligned}
f_0(x_0)-\inf\{f_0(x): x\in \Omega\}&\geq f_0(x_0)-f_0(x_k)\\
&\geq \sum^{k-1}_{i=0}\|x_{i}-x_{i+1}\|^2
\end{aligned}
\end{equation}
for any $k$, where the first inequality is due to the monotonicity of $\{f_0(x_k)\}_{k\geq 1}$ (see Theorem \ref{thm: main_theorem}) and the second one is derived through a telescoping sum. Since, due to~\eqref{eq:square_convergence}, the sum $\sum^{k}_{i=0}\|x_{i}-x_{i+1}\|^2$ converges, the increment $\|x_{k}-x_{k+1}\|^2$ goes to $0$ as $k\rightarrow \infty$. Since for $k<K$, $f_0(x_{k})-f_0(x_{k+1})\geq \mu \xi^2$, we have that $ f_0(x_0)-f_0(x_K)\geq K\mu\xi^2$. Hence,  \eqref{eq: complexity} holds.
\hfill$\blacksquare$

\section{Numerical Results}
\label{sec: experiments}
In this section, we present numerical experiments to test the performance of the SZO-QQ algorithm\footnote{All the numerical experiments have been executed on a PC with an Intel Core i9 processor.}.
\subsection{Solving an unknown convex QCQP}
To test SZO-QQ and compare it with existing zeroth-order methods, we consider the following convex {\color{black}QCQP},
\[
    \begin{aligned}
    \min_{x \in \mathbb{R}^2} & \quad  f_0(x)=0.1\times (x^{(1)})^2+x^{(2)}\\
     \text{subject to}  & \quad f_1(x) = -x^{(1)} \leq 0,\\
     &\quad f_2(x)=x^{(2)}-1\leq 0,\\
    & \quad f_3(x) = (x^{(1)})^2-x^{(2)}\leq 0.
    \end{aligned}
\]
We assume that the functions $f_i(x), i=1,2,3,$ are unknown but their values can be measured and a strictly feasible initial point $x_0 = \begin{bmatrix}0.9 & 0.9\end{bmatrix}^\top$ is given. The unique optimum $\begin{bmatrix}0&0\end{bmatrix}^\top $ is not strictly feasible. Testing on this example allows us to see whether Algorithm~\ref{alg:sslo} stays safe and whether the convergence is fast when the iterates are close to the feasible region boundary.

We run SZO-QQ and compare its result with that of s0-LBM\cite{usmanova2019log}, Extremum Seeking \cite{arnold2015model} and SafeOptSwarm~\footnote{SafeOptSwarm is a variant of SafeOpt (introduced in Section \ref{sec: introduction}). The former adds heuristics to make SafeOpt in \cite{sui2015safe} more tractable for slightly higher dimensions.}~\cite{berkenkamp2016safe}. Except Algorithm 1, s0-LBM and SafeOptSwarm are provably safe for a high probability during exploration. Only SZO-QQ and s0-LBM require Assumption \ref{ass: smoothness}. For these two methods, by trial and error (see Remark \ref{rmk: LM}), we adopt $L_{i} = 5$ and $M_{i}= 3$ for any $i\geq 1$. The penalty coefficient $\mu$ of Algorithm \ref{alg:sslo} in $\hyperlink{starstar}{(*)}$ is set to be $ 0.001$. Both s0-LBM and Extremum Seeking are log-barrier-based and we use the reformulated unconstrained problem $\min_x f_{\mathrm{log}}(x,\mu_{\log})$, where $f_{\mathrm{log}}(x,\mu_{\log}):=f_0(x) - \mu_{\mathrm{log}}\sum^4_{i=1}\log(-f_i(x))$, and $\mu_{\log} = 0.001$.
\begin{figure}[htbp!]
    \centering
    \includegraphics[width=0.5\textwidth]{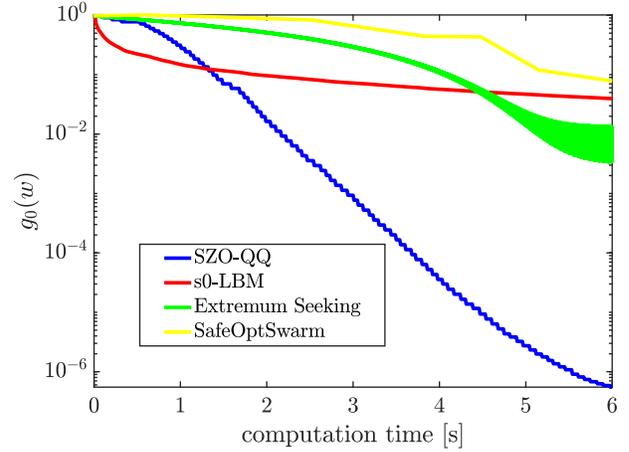}
    \caption{Objective value as a function of computation time.}
    \label{fig:objective_decrease}
\end{figure}

In Fig. \ref{fig:objective_decrease}, we show how the objective function values decrease with the computation time. During the experiments, none of the implementations violates the constraints. Regarding the convergence to the minimum, we see that s0-LBM is the fastest in the first 1.5 seconds, which is due to the low complexity of each iteration. However, SZO-QQ achieves a better solution afterwards. In the first 6 seconds, Algorithm~\ref{alg:sslo} already shows a clear trend of convergence to the optimum, which is consistent with Theorem \ref{thm: main_theorem}.

s0-LBM slows down when the iterates are close to the boundary of the feasible set (see Remark \ref{rmk: remark on the local estimation} and Example~\ref{ex: lessconservativeness}). Meanwhile, the slow convergence of Extremum Seeking is due to its small learning rate. The concept of learning rate is analogous to the step size in gradient descent. If the learning rate is large, the iterates might be brought too close to the boundary of the feasible set and then the perturbation added by this method (most visible in Fig. \ref{alg:sslo} between $5$ sec and $6$ sec) would lead to the violation of the constraints. These considerations constitute the main dilemma in parameter tuning for Extremum Seeking. In contrast, the exploration of the unknown functions in SafeOptSwarm is based on Gaussian Process (GP) regression models instead of local perturbations. Since SafeOptSwarm does not exploit the convexity of the problem, it maintains a safe set and tries to expand it for finding the global minimum. This is why many samples are close to the boundary of the feasible set even though they are far away from the optimum. This phenomenon is also observed in \cite{turchetta2019safe}. These samples along with the computational complexity of GP regression are the main reason of the slow convergence with SafeOptSwarm. 

Moreover, we also run s0-LBM and Extremum Seeking with different penalty coefficients $\mu_{\mathrm{log}}$ to check whether the slow convergence is due to improper tuning. We see that with larger $\mu_{\mathrm{log}}$ the performances of the log-barrier-based methods deteriorates. This is probably because the optimum of the unconstrained problem $f_{\log}(x,\mu_{\log})$ deviates more from the optimum $(0,0)$ as $\mu_{\mathrm{log}}$ increases. With smaller $\mu_{\mathrm{log}}$,  the Extremum Seeking method leads to constraint violation while the performance of s0-LBM barely changes. 

\subsection{Open-loop optimal control with unmodelled disturbance}
This section shows an application of SZO-QQ for solving constrained optimal control problem. We consider a non-linear system with dynamics $x_{k+1} = Ax_k+Bu_k + \delta(x_k)\begin{bmatrix}1&0\end{bmatrix}^\top$, where $x_k = [x^{(1)}_k\; x^{(2)}_k]^\top\in \mathbb{R}^2$ for $k\geq 0$, $x_0 = \begin{bmatrix} 1 & 
1\end{bmatrix}^\top,$ the matrices
$$A = \begin{bmatrix}1.1 & 1\\-0.5& 1.1 \end{bmatrix}, B  = \begin{bmatrix}1 & 0 \\0& 1 \end{bmatrix}$$
are known but the disturbance $\delta(x):=0.1*(x^{(2)})^2$ is unmodelled. We aim to design the input $u_k$ for $i=$ $0,\ldots,5,$ to minimize the cost
$\mathcal{F}(\{x_k\}_{0\leq k\leq 5 },\{u_k\}_{0\leq k\leq 5 }):= \sum^{5}_{k=0}\left(x_{k+1}^\top Q x_{k+1} +u_k^\top R u_k\right) $ where $Q= 0.5\,\mathbf I_2$ and $R = 2\,\mathbf I_2$ with identity matrix $\mathbf I_2\in\mathbb R^{2\times 2}$ while enforcing $\|x_k\|_\infty\leq 0.7$ and $\|u_k\|_\infty\leq 1.5$ for $k\geq 1$. Since we assume all the states are measured, we can evaluate the objective and constraints. This constrained optimal control problem represents a class of problems where a gap exists between the assumed and true models and safety (i.e., system trajectory verifying the constraints) is essential. These problems, called Safe Learning Control Problem (SLCP), are surveyed in \cite{brunke2022safe}. In this example, we assume to have a feasible sequence of inputs $\{u_k\}_{0\leq k \leq 5}$ (as in Assumption \ref{ass: strict feasible}) that leads to a safe trajectory and results in a cost of 6.81. Different from the settings in the model-based methods \cite{fan2020deep,hewing2019cautious} for SLCP, we don't assume that the safe trajectory is sufficient for identifying the disturbance $\delta(x)$ with small error bounds. If the error bounds are huge, the robust control problems formulated in \cite{fan2020deep,hewing2019cautious} may become infeasible.

If we only take into account the nominal model (i.e., neglect the nonlinear disturbance $\delta(x)$), the optimal control problem is a Quadratic Program. The minimal cost is 5.78 but the input sequence fails to satisfy the output constraints.

We consider applying SZO-QQ to sample only feasible trajectories and derive a feasible solution that achieves a low cost. To this aim, we utilize the known feasible input sequence as the initial solution. By setting $L_i=M_i=20$ for $i\geq 1$, $\mu = 10^{-4}$ and $\xi = 10^{-4}$, we implement Algorithm~1 to further decrease the cost resulting from the initial trajectory and derive within 57 seconds of computation an input sequence that satisfies all the constraints and achieves a cost of 5.96. This result is the same as that derived through assuming the disturbance $\delta(\cdot)$ is known and applying the nonlinear programming solver IPOPT~\cite{wachter2006implementation}, which is consistent with Theorem \ref{thm: optimality} on the convergence to a KKT point. Moreover, if we set $\xi=3*10^{-3}$, we derive, within 12 seconds, an input sequence that achieves a cost of 6.00, which shows that one may greatly reduce the computational time without ending up with a much larger cost.

Observe that the problem studied is not convex due to the disturbance. A detailed analysis of the performance of SZO-QQ for solving non-convex problems is left for future work.

\section{Conclusions}
For safe black-box convex optimization problems, we proposed a method, SZO-QQ, using samples of the objective and constraint functions to iteratively optimize over local feasible sets. Each iteration of our algorithm involves a QCQP subproblem, which can be solved efficiently. We showed that one of the accumulation points of the iterates is a minimizer. From the numerical experiments, we see our method can be faster than existing ones including s0-LBM, SafeOptSwarm and Extremum Seeking. The fast convergence of the algorithm motivates several directions of research, including the characterization of the suboptimality for the solution returned by the algorithm, the theoretical analysis of SZO-QQ when applied to non-convex problems and the study of its performance when samples of the cost function and the constraints are noisy.

\appendix
\subsection{Proof of Theorem \ref{thm: main_theorem}}
\label{sec: 1st argument of thm2}

Given any $k$, we have $x_k\in \mathcal{S}(x_k)$ and $x_{k+1} = \argmin_{x\in \mathcal{S}(x_{k})} f_0(x)+\mu \|x-x_k\|^2$. Thus, 
\begin{equation*}
\begin{aligned}
&f_0(x_{k+1})+\mu\|x_{k+1}-x_{k}\|^2 \\
\leq\;&  f_0(x_{k})+\mu\|x_{k}-x_{k}\|^2= f_0(x_k).
\end{aligned}
\end{equation*}
Then for $k\geq 0$, $f_0(x_k)\leq f_0(x_0)< \beta$ (see Assumption \ref{ass: bounded_sublevel}).

Now we know $\{x_k\}_{k\geq 1}$ is within the set $\mathcal{P}_\beta$. Due to the continuity of the functions $f_i(x)$ for $i\in \mathbb{Z}^m_0$ and the boundedness of the set $\mathcal{P}_\beta$, the set $\mathcal{P}_\beta$ is compact. Hence, $\{x_k\}_{k\geq 1}$ has at least an accumulation point $x_c$.

Since $f_0(x)$ is continuous and $\Omega$ is closed, $\inf_{x\in \Omega}f_0(x)$ $>-\infty$ and the non-increasing sequence $(f_0(x_k))_{k\geq 1}$ converges to $\inf_{k}f_0(x_k)$.  Considering that $x_c$ is the limit point of a subsequence, we have $f_0(x_c)\geq \inf_{k}f_0(x_k)$. Due to the fact that, for any $k$, $f_0(x_c)\leq f_0(x_k)$, we have $f_0(x_c)\leq \inf_{k}f_0(x_k)$. Thus $\{f_0(x_k)\}_{k\geq 1}$ converges to $f_0(x_c)$.

\subsection{Statement and proof of Lemma \ref{prop: active set non-empty}}
\label{sec: state_proof_lemma_nonempty_active_set}
\begin{lemma}
\label{prop: active set non-empty}
If $x_c$ is an accumulation point of the sequence $\{x_k\}_{k\geq 1}$ and $\nabla f_0(x_c)\neq 0$, the set $\mathcal{A}(x_c)$ is nonempty.
\end{lemma}
\textbf{Proof.}
We prove the lemma by contradiction. We assume $f_i(x_c)<0$ for all $i$ and let $(x_{k_j})_j$ be the subsequence converging to $x_c$.

The first claim to prove is that for some $r_c>0$ there exists $j_c$ such that, for any $j>j_c$, $\mathcal{B}_{r_c}(x_c)\subset \mathcal{S}(x_{k_j})$. To see this, we let $j_c$ be such that for any $j\geq j_c$ and any $i$,
\begin{equation}
\label{eq: large_enough_j}
\begin{aligned}
f_i(x_{k_j})&<\frac{1}{2}f_i(x_c),\\
2L_{\max}\|x_c-x_{k_j}\|+4M_{\max}\|x_c-x_{k_j}\|^2&<-\frac{1}{4}f_i(x_c).
\end{aligned}    
\end{equation}

Since $\nu_{k_j}$ converges to 0, ${\nabla}^{\nu^*_{k_j}} f_i\left(x_{k_j}\right)$ converges to ${\nabla} f_i\left(x_{k_j}\right)$. Thus, without loss of generality, we can assume, for any $j\geq j_c$, $\|{\nabla}^{\nu^*_{k_j}} f_i\left(x_{k_j}\right)\|\leq 2L_{\max}$. We also let $r_c>0$ be such that 
$$2L_{\max}r_c+4M_{\max}r_c^2<-\frac{1}{4}f_i(x_c).$$
Then, for any $x'_c\in \mathcal B_{r_c}(x_c)$ and $j\geq j_c$ we have 
\[
\begin{aligned}
& f_i(x_{k_j})+{\nabla}^{\nu^*_{k_j}} f_i\left(x_{k_j}\right)(x'_c-x_{k_j})+2M_i\|x'_c-x_{k_j}\|^2\\
< & -\frac{1}{2}f_i(x_c) + 2L_{\max}\|x'_c-x_{k_j}\|+2M_{\max}\|x'_c-x_{k_j}\|^2\\
< & -\frac{1}{2}f_i(x_c) + 2L_{\max}\|x_c-x_{k_j}\|+ 4M_{\max}\|x_c-x_{k_j}\|^2+\\
 & \text{ }2L_{\max}\|x_c-x'_c\| + 4M_{\max}\|x_c-x'_c\|^2\\
 < & -\frac{1}{2}f_i(x_c) + \frac{1}{4}f_i(x_c) + \frac{1}{4}f_i(x_c) = 0,
\end{aligned}
\]
where the second inequality is due to the triangle inequality while the third is derived using \eqref{eq: large_enough_j}. This means for any $j\geq j_c$  we have $\mathcal{B}_{r_c}(x_c)\subset S(x_{k_j})$.

The second claim is that in $\mathcal{B}_{r_c}(x_c)$ we can find $x$ such that $f_0(x)+ \mu\|x-x_c\|^2<f_0(x_c)$. To see this we let $x = x_c-t\nabla f_0(x_c)$, where $t\in \mathbb{R}_{\geq 0}$ is a parameter. From \eqref{eq: meanvaluetheorem}, we have that
\[
\begin{aligned}
f_0(x) & \leq f_0(x_c)+\nabla f_0(x_c)^\top(x-x_c)+M_0\|x-x_c\|^2\\
&= f_0(x_c) - t\|\nabla{f}_0(x_c)\|^2+t^2M_0\|\nabla{f}_0(x_c)\|^2.
\end{aligned}
\]
Hence, for $t<\min \{(M_0+\mu)^{-1},r_c/|\nabla f_0(x_c )|\}$, we have 
\[
f_0(x)+ \mu\|x-x_c\|^2<f_0(x_c).
\]
Therefore, there exists $j_c\geq 1$ such that, for any $j\geq j_c$, $f_0(x)+ \mu\|x-x_{k_j}\|^2< f_0(x_{c}).$ Since $x\in S(x_{k_j})$, we have
\begin{equation*}
\begin{aligned}
&f_0(x_{k_j+1})+ \mu\|x_{k_j+1}-x_{k_j}\|^2\\\leq&  f_0(x)+ \mu\|x-x_{k_j}\|^2\,<\, f_0(x_{c}).
\end{aligned}
\end{equation*}
which leads to the fact that $f_0(x_{k_j+1})<f_0(x_c)$ and $$\underset{k\rightarrow\infty}{\lim}f_0(x_k)<f_0(x_c),$$
This result contradicts Theorem \ref{thm: main_theorem}. Thus, $\mathcal{A}(x_c)$ is not empty.\hfill$\blacksquare$

\subsection{Statement and proof of Lemma \ref{claim: more_optimal_point}}
\label{sec: state_proof_more_optimal_point}
\begin{lemma}
\label{claim: more_optimal_point}
If $x_c$ is not an optimum, $\nabla f_0(x_c)\neq 0$ and Assumption \ref{ass: constraint qualification} is verified, the following properties hold:
\begin{itemize}
    \item There exists $x_s\in\Omega$ such that $f_0(x_s)<f_0(x_c)$ and $\nabla f_i(x_c)^\top (x_s-x_c)<0$, $\forall i\in \mathcal{A}(x_c)$;
    \item Define $x(t):=x_c+t(x_s-x_c), t\in [0,1]$ where $x_s$ was given above. Then, there exists $\bar{t}\leq 1$ such that, for $0<t<\bar{t}$, $f_0(x(t))+\mu \|x(t)-x_c\|^2 <f_0(x_c)$;
    \item There exists $0<t_p<\bar{t}$ such that, for any $\tau>0$ and $k_0\in\mathbb{Z}_{\geq 1}$, we can find $k_p> k_0$ satisfying $\|x_{k_p}-x_c\|\leq\tau$ and $ x(t_p)\in \mathcal{S}(x_{k_p})$.
\end{itemize} 
\end{lemma}
\textbf{Proof.} Convexity of $\Omega$ along with Assumption \ref{ass: strict feasible} on strict feasibility leads to the following property.
\begin{property}
\label{fct: interior_point}
For any $x\in \Omega$ and any $r>0$, there exists an interior point of $\Omega$ in $\mathcal B_r(x)$, i.e., there exists $x'\in \mathcal B_r(x)$ and $r'>0$ such that $\mathcal B_{r'}(x')\subset \Omega$.
\end{property}

If $x_c$ is not a minimizer, we let $x_s\in\Omega$ be such that $f_0(x_s)<f_0(x_c)$.
According to Property~\ref{fct: interior_point} and the assumption that $\nabla f_i(x_c)\neq 0$ for $i\in \mathcal{A}(x_c)$, we can find an interior point $x'_s$ in the neighborhood of $x_s$ such that $f_0(x'_s)<f_0(x_c)$ and another point $x''_s\in \Omega$ in the neighborhood of $x'_s$ such that $f_0(x''_s)<f_0(x_c)$ and $-\nabla f_i(x_c)^\top (x''_s-x_c)\neq 0, \forall i\in \mathcal{A}(x_c)$. Since $x(t):=tx_c+(1-t)x''_s\in \Omega$ for $0<t<1$, there exists $g(t)$ such that 
\begin{equation*}
\begin{aligned}
\lim_{t\rightarrow \infty} \frac{g(t)}{t} &= 0,\\
f_i(x(t)) &= f_i(x_c)+\nabla f_i(x_c)^\top (x''_s-x_c)t + g(t)\leq 0.
\end{aligned}
\end{equation*}
Therefore, $-\nabla f_i(x_c)^\top (x''_s-x_c)>0, \forall i\in \mathcal{A}(x_c)$. 
\vspace{0.1cm}

\begin{spacing}{1.05}
Due to convexity, $f_0(x_s)\geq f_0(x_c) + \nabla f_0(x_c)^\top(x_s-x_c)$, which leads to the fact that $\nabla f_0(x_c)^\top(x_s-x_c)<0$ and there exists a constant $\sigma <0$ such that $\nabla f_0(x_c)^\top (x(t)-x_c) = t\sigma $. Since $f_0(x(t))\leq f_0(x_c)+t\sigma+t^2M_0\|x_s-x_c\|^2$, we let $\bar{t} = -\sigma/((M_0+\mu)\|x_s-x_c\|^2)$ and have, for $0<t<\bar{t}$, $f_0(x(t))+\mu \|x(t)-x_c\|^2 <f_0(x_c)$.
\end{spacing} 

To check the feasibility of $x(t)$ with respect to $\mathcal{S}_i(x_k)$ for $i\in \mathcal{A}(x_c)$, we first notice that $$\mathcal{S}_i(x_k) = \mathcal{B}_{r^k_i(x_k)}(O^k_i(x_k)), \text{where for a strictly feasible }x, $$ $$O^k_i(x) := x-\frac{\nabla f^{\nu^*_{k}}_i(x)}{2M_i} \text{, } \left(r^k_i(x)\right)^2 :=-\frac{f_i(x)}{M_i}+\frac{\|\nabla f_i^{\nu^*_{k}}(x)\|^2}{4M_i^2}.$$
Furthermore, for $x_c$, we define 
$$O_i(x_c) := x_c-\frac{\nabla f_i(x_c)}{2M_i} \text{, } \left(r_i(x_c)\right)^2 :=-\frac{f_i(x_c )}{M_i}+\frac{\|\nabla f_i(x_c)\|^2}{4M_i^2}.$$ With $\bar{\bar{t}}_i:=\argmin_{t\in \mathbb{R}^d}\|O_i(x_c)-x(t)\|$, by considering that
\begin{equation*}
\|O_i(x_c)-x_c+x_c-x(\bar{\bar{t}}_i)\|  
=\|\frac{\nabla f_i(x_c)}{2M_i} + \bar{\bar{t}}_i(x_s-x_c)\|
\end{equation*}
and the fact that $\|O_i(x_c)-x(\bar{\bar{t}}_i)\|\leq r_i(x_c)$, we know 
$$(r_i(x_c))^2+\|\bar{\bar{t}}_i(x_s-x_c)\|^2+\frac{\bar{\bar{t}}_i\nabla f_i(x_c)^\top(x_s-x_c)}{2M_i}\leq (r_i(x_c))^2.$$
Since $-\nabla f_i(x_c)(x_s-x_c)>0$, we have $\bar{\bar{t}}_i>0$. Then for any $0<\tilde{t}_i\leq \min\{\bar{\bar{t}}_i,\bar{t},1\}$, we have $\|O_i(x_c)-x(\tilde{t}_i)\|<r_i(x_c)$ and thus $x(\tilde{t}_i)$ is an interior point in $ \mathcal{B}_{r_i(x_c)}(O_i(x_c)), \forall i\in \mathcal{A}(x_c)$. Consequently, we can find $t_p \leq \bar{t}$ such that $x(t_p)$ is an interior point in the ball $ \mathcal{B}_{r_i(x_c)}(O_i(x_c)),\forall i\in \mathcal{A}(x_c),$ and $f_0(x(t_p))<f_0(x_c)$.

We let $\{x_{k_j}\}_{j\geq 1}$ be the subsequence of $\{x_k\}_{k\geq 1}$ that converges to $x_c$. Since $\nu^*_{k_j}$ converges to 0, $ O^{k_j}_i(x_{k_j})$ converges to $O_i(x_c)$ and $r^{k_j}_i(x_{k_j})$ converges to $r_i(x_c)$ for any $i\in\mathcal{A}(x_c)$. Because $x(t_p)$ is an interior point of $ \mathcal{B}_{r_i(x_c)}(O_i(x_c)),\forall i\in \mathcal{A}(x_c),$ we know that for any $j_0$ there exists $j>j_0$ such that $x_{k_j}$ is arbitrarily close to $x_c$ and $x(t_p)\in \mathcal{S}(x_{k_j})$. Thus, for any $k_0$ and $\tau$ there exists $k_p>k_0$ such that $\|x_{k_p}-x_c\|\leq \tau$ and $x(t_p)\in \mathcal{S}(x_{k_p})$.\hfill$\blacksquare$

\subsection{Proof of Theorem \ref{thm: optimality}}
\label{sec: proof_of_optimality}
We let $x_c$ be an accumulation point of $\{x_k\}_{k\geq 1}$ such that $\nabla f_i(x_c)\neq 0,\forall i\in \mathcal{A}(x_c)$ (see Assumption \ref{ass: constraint qualification}). If $\nabla f_0(x_c) = 0$, then with $\lambda_c = 0$, the pair $(x_c,\lambda_c)$ is a KKT point and thus $x_c$ is minimizer. In the following, we assume $\nabla f_0(x_c) \neq 0$.

From Lemma \ref{prop: active set non-empty}, we know
\begin{equation*}
\label{eq: conditions_for_contradiction}
\nabla f_0(x_c)\neq 0, \; \mathcal{A}(x_c)\neq \emptyset, \; \text{and }  \nabla f_i(x_c)\neq 0 \text{ for } i\in\mathcal{A}(x_c).
\end{equation*}
Then, we prove the optimality of $x_c$ by contradiction. By assuming $x_c$ is not an optimum, we know the three points in Lemma \ref{claim: more_optimal_point} hold. Therefore, for any $\tau>0$ and $k_0\in\mathbb{Z}_{\geq 1}$ there exist $k_p>k_0$ and $x'\in \mathcal{S}(x_{k_p})$ verifying $\|x_{k_p}-x_c\|\leq\tau$ and $f_0(x')+\mu \|x'-x_c\|^2 <f_0(x_c)$. By letting $\tau$ be sufficiently small, we have $k_p$ and $x'\in \mathcal{S}(x_{k_p})$ verifying 
\[
\mu\|x_c-x_{k_p}\|^2<\frac{1}{2}\left(f_0(x_c)-f_0(x(t_p))-\mu \|x(t_p)-x_c\|^2\right),
\] 
and $f_0(x')+\mu \|x'-x_c\|^2 <f_0(x_c).$
Therefore, we have 
\begin{equation*}
\begin{aligned}
& f_0(x')+\mu \|x'-x_{k_p}\|^2  \\
\leq & f_0(x')+\mu \|x'-x_{c}\|^2 + \mu \|x_{c}-x_{k_p}\|^2\\
<&\frac{1}{2}f_0(x_c)+\frac{1}{2}\left(f_0(x')+\mu \|x(t_p)-x_c\|^2\right)\\
<& f_0(x_c).
\end{aligned}
\end{equation*}
Therefore, we have 
\begin{equation*}
\begin{aligned}
&f_0(x_{k_p+1})+ \mu\|x_{k_p+1}-x_{k_p}\|^2\\
\leq &  f_0(x')+\mu \|x'-x_{k_p}\|^2\;<\; f_0(x_{c}),
\end{aligned}
\end{equation*}
which contradicts the monotonicity properties shown in Theorem \ref{thm: main_theorem}. Now, we have shown that $x_c$ is optimum.

For any accumulation point $x_c$, no matter whether $\nabla f_i(x_c)\neq 0,\forall i\in \mathcal{A}(x_c)$ holds, we know $f_0(x_c) = \lim_{k\rightarrow \infty}f_0(x_k)$ due to Theorem \ref{thm: main_theorem}. Thus, any accumulation point $x_c$ is a minimizer of \eqref{eq: optimization problem}.

\bibliographystyle{ieeetr}
\bibliography{library} 
\end{document}